\documentclass[a4paper,notitlepage,11pt]{article}%
\usepackage{amssymb}
\usepackage{graphicx}
\usepackage{amsmath}
\usepackage{amsthm}
\usepackage{amsfonts}%
\providecommand{\U}[1]{\protect\rule{.1in}{.1in}}
\theoremstyle{plain}
\newtheorem{theorem}{Theorem}[section]

\newtheorem{corollary}[theorem]{Corollary}

\newtheorem{lemma}[theorem]{Lemma}

\theoremstyle{definition}

\newtheorem{remark}[theorem]{Remark}
\numberwithin{equation}{section}
\numberwithin{theorem}{section}
\ifx\pdfoutput\relax\let\pdfoutput=\undefined\fi
\newcount\msipdfoutput
\ifx\pdfoutput\undefined\else
\ifcase\pdfoutput\else
\ifx\paperwidth\undefined\else
\ifdim\paperheight=0pt\relax\else\pdfpageheight\paperheight\fi
\ifdim\paperwidth=0pt\relax\else\pdfpagewidth\paperwidth\fi
\fi\fi\fi
\begin{document}

\title{One-dimensional singular problems involving the $p$-Laplacian and
nonlinearities indefinite in sign\thanks{2010 \textit{Mathematics Subject
Clasification}. 34B16; 34B18, 34B15, 34C25.}  \thanks{Partially supported by Secyt-UNC and
CONICET.} }
\author{U. Kaufmann, I. Medri\thanks{\textit{E-mail addresses. }%
kaufmann@mate.uncor.edu (U. Kaufmann, Corresponding Author),
medri@mate.uncor.edu (I. Medri).}
\and \noindent\\{\small FaMAF, Universidad Nacional de C\'{o}rdoba, (5000) C\'{o}rdoba,
Argentina}}
\maketitle

\begin{abstract}
Let $\Omega$ be a bounded open interval, let $p>1$ and $\gamma>0$, and let
$m:\Omega\rightarrow\mathbb{R}$ be a function that may change sign in $\Omega
$. In this article we study the existence and nonexistence of positive
solutions for one-dimensional singular problems of the form $-(\left\vert
u^{\prime}\right\vert ^{p-2}u^{\prime})^{\prime}=m\left(  x\right)
u^{-\gamma}$ in $\Omega$, $u=0$ on $\partial\Omega$. As a consequence we also
derive existence results for other related nonlinearities.

\end{abstract}

\section{Introduction}

For $a<b$, let $\Omega:=(a,b)$, and let $\gamma>0$. Let $p\in\left(
1,\infty\right)  $ and $m\in L^{p^{\prime}}\left(  \Omega\right)  $ (where as
usual we define $p^{\prime}$ by $1/p+1/p^{\prime}=1$) be a possibly sign
changing function, and consider the problem%

\begin{equation}
\left\{
\begin{array}
[c]{ll}%
-\left(  \left\vert u^{\prime}\right\vert ^{p-2}u^{\prime}\right)  ^{\prime
}=m\left(  x\right)  u^{-\gamma} & \text{in }\Omega\\
u>0 & \text{in }\Omega\\
u=0 & \text{on }\partial\Omega.
\end{array}
\right.  \label{prob}%
\end{equation}
One-dimensional singular problems involving the $p$-Laplacian like
(\ref{prob}) arise in applications such as non-Newtonian fluid theory or the
turbulent flow of a gas in a porous medium (cf. \cite{siam}, \cite{gh4}), and
they have been widely studied over the years if $m$ is nonnegative. We cite,
among many others, the papers \cite{wang}, \cite{lu}, \cite{agar},
\cite{agar2003}, \cite{amc}, \cite{ma}, \cite{sun}. However, to the best of
our knowledge, there are no results available in the literature when $m$ is
allowed to change sign in $\Omega$. Let us note that if $m$ has an indefinite
sign (\ref{prob}) becomes a much more involved problem. In fact, (when $m$
changes sign) these problems are quite intriguing even when (\ref{prob}) is
sublinear (i.e., $\gamma\in\left(  1-p,0\right)  $), and only lately existence
of positive solutions have been obtained in this case (see \cite{plap} for
$p\in\left(  1,\infty\right)  $, and \cite{ejde} and its references for the
special case $p=2$).

On the other side, for the Laplace operator (that is, $p=2$) the problem
(\ref{prob}) has been recently considered in \cite{jmaa} for sign changing
$m$'s. Our aim in this article is to establish similar results in the general
situation $1<p<\infty$, adapting and extending the approach developed in
\cite{jmaa} combined also with some of the ideas in \cite{plap}. Let us
mention that this is far from being trivial due to the nonlinearity of the
$p$-Laplacian and its corresponding solution operator. Moreover, we remark
that some of the conditions presented in this paper improve the ones found in
\cite{jmaa} for the laplacian operator.

In order to derive our results we shall mainly rely on the well-known sub and
supersolution method. The major difficulty here (as with various nonlinear
problems with indefinite nonlinearities) is to find a (strictly) positive
subsolution. We shall provide such subsolution by means of Schauder's fixed
point theorem applied to some related nonlinear problems. More precisely, in
Theorem \ref{nuevo} (i) we shall give a sufficient condition on $m$ that
assures the existence of solutions of (\ref{prob}) for all $\gamma>0$ small
enough, and further conditions are stated in Theorem \ref{nuevo} (ii) without
the smallness restriction on $\gamma$ (see also Remark \ref{dos} below).

On the other hand, two necessary conditions on $m$ are exhibited in Theorem
\ref{necesaria} (see also Remark \ref{rito}). Let us point out that the first
of the aforementioned sufficient conditions on $m$ turns out to be also
\textquotedblleft almost\textquotedblright\ necessary (compare (\ref{nece})
with (\ref{ne}), and see the last paragraph in Remark \ref{rito}). Finally, as
a consequence of the above theorems, we shall prove in Corollary \ref{fff} an
existence result for singular nonlinearities of the form $m\left(  x\right)
f\left(  u\right)  $ with no monotonicity nor convexity assumptions on $f$.

We conclude this introduction with some few comments on some related open
interesting problems. Based on the results in \cite{plap} for the analogous
sublinear problem, we think that similar theorems to the ones proved here
should still be true replacing the $p$-Laplacian by operators of the form
\[
\mathcal{L}u=-\left(  \left\vert u^{\prime}\right\vert ^{p-2}u^{\prime
}\right)  ^{\prime}+c\left(  x\right)  \left\vert u\right\vert ^{p-2}u,
\]
where $c\geq0$ in $\Omega$. We note however that, for instance, the proof of
the key Lemma \ref{plap} does not work in this case and it is not clear how to
adapt it. Also somehow similar results should be valid for the analogous
$n$-dimensional problem (in fact, this occurs when $p=2$ (and $c\equiv0$), see
\cite{jmaa}, Section 4; and also \cite{gh1}, \cite{gh2}, \cite{gh3} for
related elliptic problems), and in our opinion proving this if $p\not =2$ is
not a trivial task. Let us finally mention that in the one-dimensional case
one could also consider (\ref{prob}) with the so-called $\phi$-Laplacian in
place of the $p$-Laplacian, that is, taking $\left(  \phi\left(  u^{\prime
}\right)  \right)  ^{\prime}$ instead of the $p$-Laplacian, where $\phi$ is an
increasing odd homeomorphism with $\phi\left(  \mathbb{R}\right)  =\mathbb{R}$
(for singular problems with the $\phi$-Laplacian we refer to the book
\cite{fi}, Part II).

\section{Preliminaries}

For $1<p<\infty$, let $\mathcal{L}$ be the differential operator given by
\[
\mathcal{L}v:=-\left(  \left\vert v^{\prime}\right\vert ^{p-2}v^{\prime
}\right)  ^{\prime}.
\]
We start collecting some necessary facts concerning the problem%
\begin{equation}
\left\{
\begin{array}
[c]{ll}%
\mathcal{L}v=h\left(  x\right)  & \text{in }\Omega\\
v=0 & \text{on }\partial\Omega.
\end{array}
\right.  \label{g}%
\end{equation}

\begin{remark}
\label{uno} Let $h\in L^{q}(\Omega)$, $q>1$. It is well known that (\ref{g})
admits a unique solution $v\in C^{1}\left(  \overline{\Omega}\right)  $ such
that $\left\vert v^{\prime}\right\vert ^{p-2}v^{\prime}$ is absolutely
continuous and that the equation holds in the pointwise sense. In fact, if%
\[
\varphi_{p}\left(  t\right)  :=\left\vert t\right\vert ^{p-2}t\qquad\text{for
}t\not =0,\qquad\varphi_{p}\left(  0\right)  :=0,
\]
and $\varphi_{p}^{-1}$ denotes its inverse, it can be seen that%
\begin{equation}
v\left(  x\right)  =\int_{a}^{x}\varphi_{p}^{-1}\left(  c_{h}-\int_{a}%
^{y}h\left(  t\right)  dt\right)  dy, \label{v1}%
\end{equation}
where $c_{h}$ is the unique constant such that $v\left(  b\right)  =0$ (see
e.g. \cite{mana}, Section 2). Furthermore, the solution operator $\mathcal{S}$
satisfies that $\mathcal{S}:L^{q}(\Omega)\rightarrow C^{1}(\overline{\Omega})$
is continuous (e.g. Lemma 2.1 \ in \cite{ma1} or Lemma 4.2 in \cite{ma2}) and
$\mathcal{S}:L^{q}(\Omega)\rightarrow C(\overline{\Omega})$ is compact (cf.
\cite{mana}, Corollary 2.3). $\blacksquare$
\end{remark}

The so-called weak comparison principle shall be repeatedly used along the
paper, and so we state it here for the reader's convenience (for a proof, see
for instance \cite{les}, Corollary 6.5.3).

\begin{lemma}
\label{w}Let $u,v\in W_{0}^{1,p}\left(  \Omega\right)  $ be such that $u\leq
v$ on $\partial\Omega$ and $\mathcal{L}u\leq$ $\mathcal{L}v$ in weak sense in
$\Omega$, that is,
\[
\int_{a}^{b}\left\vert u^{\prime}\right\vert ^{p-2}u^{\prime}\varphi^{\prime
}\leq\int_{a}^{b}\left\vert v^{\prime}\right\vert ^{p-2}v^{\prime}%
\varphi^{\prime}\qquad\text{for all }0\leq\varphi\in W_{0}^{1,p}\left(
\Omega\right)  .
\]
Then $u\leq v$ in $\Omega$.
\end{lemma}

The next remark compiles some properties concerning the first eigenvalue of
the $p$-Laplacian and its corresponding eigenfunctions.

\begin{remark}
\label{autos}There exists a first eigenvalue $\lambda_{1}\left(
\Omega\right)  >0$ and $\Phi\in W_{0}^{1,p}\left(  \Omega\right)  $,
$\left\Vert \Phi\right\Vert _{L^{\infty}\left(  \Omega\right)  }=1$,
satisfying
\begin{equation}
\left\{
\begin{array}
[c]{ll}%
\mathcal{L}\Phi=\lambda_{1}\left(  \Omega\right)  \Phi^{p-1} & \text{in
}\Omega\\
\Phi>0 & \text{in }\Omega\\
\Phi=0 & \text{on }\partial\Omega.
\end{array}
\right.  \label{auto}%
\end{equation}
Moreover,
\[
\lambda_{1}\left(  \Omega\right)  =\left(  \frac{\pi_{p}}{b-a}\right)
^{p},\qquad\text{where\qquad}\pi_{p}:=\frac{2\pi\left(  p-1\right)  ^{1/p}%
}{p\sin\left(  \pi/p\right)  }\text{;}%
\]
and $\Phi$ is a multiple of the function $\sin_{p}\left(  \pi_{p}\left(
x-a\right)  /\left(  b-a\right)  \right)  $ which is strictly positive and
symmetric in $\Omega$ and increasing in $\left(  a,\left(  a+b\right)
/2\right)  $ (see e.g. \cite{les}, Section 6.3; and for the precise definition
and further properties of $\sin_{p}$, see e.g. \cite{lin} and \cite{bie},
Section 2). $\blacksquare$
\end{remark}

In the following lemma we establish some useful upper and lower bounds for
$\mathcal{S}\left(  h\right)  $. We write as usual $h=h^{+}-h^{-}$ with
$h^{+}:=\max\left(  h,0\right)  $ and $h^{-}:=\max\left(  -h,0\right)  $. We
also set
\[
\delta_{\Omega}\left(  x\right)  :=dist\left(  x,\partial\Omega\right)
=\min\left(  x-a,b-x\right)  \text{.}%
\]

\begin{lemma}
\label{plap}Let $p\in\left(  1,\infty\right)  $ and $h\in L^{q}(\Omega)$ for
some $q>1$. \newline(i) If $h\geq0$, then in $\overline{\Omega}$ it holds that%
\begin{equation}
\mathcal{S}\left(  h\right)  \leq\left(  \int_{a}^{b}h\right)  ^{1/\left(
p-1\right)  }\delta_{\Omega}. \label{lap1}%
\end{equation}
(ii) Let $I:=\left(  x_{0},x_{1}\right)  \subseteq\Omega$ and let
$x_{I}:=\left(  x_{0}+x_{1}\right)  /2$. If%
\begin{equation}
\inf_{I}h>\lambda_{1}\left(  I\right)  \max\left(  \left(  x_{I}-a\right)
^{p-1}\int_{a}^{x_{0}}h^{-},\left(  b-x_{I}\right)  ^{p-1}\int_{x_{1}}%
^{b}h^{-}\right)  , \label{hipo}%
\end{equation}
then in $\overline{\Omega}$ it holds that%
\begin{gather}
\mathcal{S}\left(  h\right)  \geq\min\left(  \mathcal{H}_{a},\mathcal{H}%
_{b}\right)  ^{1/\left(  p-1\right)  }\delta_{\Omega},\qquad\text{where}%
\label{lap2}\\
\mathcal{H}_{a}:=\frac{\inf_{I}h}{\lambda_{1}\left(  I\right)  \left(
x_{I}-a\right)  ^{p-1}}-\int_{a}^{x_{0}}h^{-},\qquad\mathcal{H}_{b}%
:=\frac{\inf_{I}h}{\lambda_{1}\left(  I\right)  \left(  b-x_{I}\right)
^{p-1}}-\int_{x_{1}}^{b}h^{-}.\nonumber
\end{gather}

\end{lemma}

\textit{Proof}. Let us prove (i). We assume here without loss of generality
that $h\not \equiv 0$. Then by the strong maximum principle (e.g.
\cite{garcia}, Theorem 2) $\mathcal{S}\left(  h\right)  >0$ in $\Omega$. We
observe now that $\varphi_{p}^{-1}=t^{1/\left(  p-1\right)  }$ for $t\geq0$
and $\varphi_{p}^{-1}=-\left\vert t\right\vert ^{1/\left(  p-1\right)  }$ if
$t<0$, and so using (\ref{v1}) we discover that $\mathcal{S}\left(  h\right)
^{\prime}\left(  x\right)  =\varphi_{p}^{-1}\left(  c_{h}-\int_{a}^{x}h\left(
t\right)  dt\right)  $ is nonincreasing because $h\geq0$. Hence,
$\mathcal{S}\left(  h\right)  $ is concave in $\Omega$ and thus it must hold
that $\mathcal{S}\left(  h\right)  ^{\prime}\left(  b\right)  <0<\mathcal{S}%
\left(  h\right)  ^{\prime}\left(  a\right)  $ and therefore
\begin{equation}
0<c_{h}<\int_{a}^{b}h\left(  t\right)  dt. \label{ch}%
\end{equation}
Noticing that $\varphi_{p}^{-1}$ is increasing and (\ref{ch}) we get that
$\mathcal{S}\left(  h\right)  ^{\prime}\left(  a\right)  ,\left\vert
\mathcal{S}\left(  h\right)  ^{\prime}\left(  b\right)  \right\vert
\leq\left(  \int_{a}^{b}h\right)  ^{1/\left(  p-1\right)  }$ and then from the
concavity of $\mathcal{S}\left(  h\right)  $ we derive (\ref{lap1}).

On the other side, let $I:=\left(  x_{0},x_{1}\right)  \subseteq\Omega$ , and
let $\lambda_{1}(I)>0$ and $\Phi>0$ with $\left\Vert \Phi\right\Vert
_{L^{\infty}\left(  I\right)  }=1$ be the corresponding normalized positive
eigenfunction for the $p$-Laplacian in $I$, that is, satisfying (\ref{auto})
with $I$ in place of $\Omega$. Suppose that (\ref{hipo}) holds (in particular,
$\inf_{I}h>0$) and fix $\lambda^{\ast}:=\lambda_{1}\left(  I\right)  /\inf
_{I}h$. In order to prove (ii) we start building some $0<u\in W_{0}%
^{1,p}\left(  \Omega\right)  $ such that $\mathcal{L}u\leq\lambda^{\ast
}h\left(  x\right)  $ in weak sense in $\Omega$. Its construction is inspired
in some of the computations made in the proofs of Theorems 3.1 and 3.5 in
\cite{plap} and \cite{ejde} respectively. Let us first point out that since
$0<\Phi\leq1$,%
\begin{equation}
\mathcal{L}\Phi=\lambda_{1}(I)\Phi^{p-1}\leq\lambda^{\ast}h\left(  x\right)
\text{\qquad in }I.\label{estre}%
\end{equation}

On the other hand, define
\begin{gather*}
c_{a}:=\frac{1}{\left(  x_{I}-a\right)  ^{p-1}}-\lambda^{\ast}\int_{a}^{x_{0}%
}h^{-},\\
v\left(  x\right)  :=\int_{a}^{x}\left(  c_{a}+\lambda^{\ast}\int_{a}^{y}%
h^{-}\right)  ^{1/\left(  p-1\right)  }dy,\text{\qquad}x\in\left[
a,x_{I}\right]  .
\end{gather*}
(Recall that $x_{I}:=\left(  x_{0}+x_{1}\right)  /2$, and note that $c_{a}>0$
due to (\ref{hipo}).) It is easy to check that $v$ is increasing and convex,
$v\left(  a\right)  =0$ and $\mathcal{L}v=-\lambda^{\ast}h^{-}\left(
x\right)  \leq\lambda^{\ast}h\left(  x\right)  $ in $\left(  a,x_{I}\right)
$. Also, (\ref{hipo}) implies that $h>0$ in $I$ and thus
\[
\left\Vert v\right\Vert _{L^{\infty}\left(  a,x_{I}\right)  }\leq\int
_{a}^{x_{I}}\left(  c_{a}+\lambda^{\ast}\int_{a}^{x_{0}}h^{-}\right)
^{1/\left(  p-1\right)  }dy=1.
\]
Similarly, if for $x\in\left[  x_{I},b\right]  $ we set
\begin{gather*}
c_{b}:=\frac{1}{\left(  b-x_{I}\right)  ^{p-1}}-\lambda^{\ast}\int_{x_{1}}%
^{b}h^{-},\\
w\left(  x\right)  :=\int_{x}^{b}\left(  c_{b}+\lambda^{\ast}\int_{y}^{b}%
h^{-}\right)  ^{1/\left(  p-1\right)  }dy,
\end{gather*}
then $w$ is decreasing and convex, $w\left(  b\right)  =0$, $\mathcal{L}%
w\leq\lambda^{\ast}h\left(  x\right)  $ in $\left(  x_{I},b\right)  $ and
$\left\Vert w\right\Vert _{L^{\infty}\left(  x_{I},b\right)  }\leq1$.

Now, since $v\left(  a\right)  =w\left(  b\right)  =\Phi\left(  x_{0}\right)
=\Phi\left(  x_{1}\right)  =0$ and $\left\Vert v\right\Vert _{\infty
},\left\Vert w\right\Vert _{\infty}\leq1=\left\Vert \Phi\right\Vert _{\infty}%
$, and since $\Phi$ is increasing in $\left[  x_{0},x_{I}\right]  $ and
decreasing in $\left[  x_{I},x_{1}\right]  $ (see Remark \ref{autos}),
reasoning as in the proof of Theorem 3.1 (i) in \cite{ejde} we find some
$\underline{x}_{0}\in\left(  x_{0},x_{I}\right)  $ and $\overline{x}_{1}%
\in\left(  x_{I},x_{1}\right)  $ such that
\begin{gather}
v(\underline{x}_{0})=\Phi(\underline{x}_{0}),\text{\qquad}\Phi\left(
\overline{x}_{1}\right)  =w\left(  \overline{x}_{1}\right)  ,\label{ollo}\\
v^{\prime}(\underline{x}_{0})\leq\Phi^{\prime}(\underline{x}_{0}%
),\text{\qquad}\Phi^{\prime}(\overline{x}_{1})\leq w^{\prime}(\overline{x}%
_{1}).\nonumber
\end{gather}
Let us define a function $u$ by $u:=v$ in $\left[  a,\underline{x}_{0}\right]
$, $u:=\Phi$ in $\left[  \underline{x}_{0},\overline{x}_{1}\right]  $ and
$u:=w$ in $\left[  \overline{x}_{1},b\right]  $. (We mention that if $x_{0}%
=a$, in order to build $u$ we only use $\Phi$ and $w$, if $x_{1}=b$ then we do
not need $w$, and if $I=\Omega$ we simply put $u=\Phi$.) Taking into account
the above paragraph, (\ref{estre}) and (\ref{ollo}), a simple integration by
parts gives that $\mathcal{L}u\leq\lambda^{\ast}h\left(  x\right)  $ in weak
sense in $\Omega$. Moreover, since $v^{\prime}\left(  a\right)  =c_{a}%
^{1/\left(  p-1\right)  }$ and $-w^{\prime}\left(  b\right)  =c_{b}^{1/\left(
p-1\right)  }$, by the convexity of $v$ and $w$ and the aforementioned
monotonicity properties of $\Phi$ it follows that
\[
u\geq\min\left(  c_{a},c_{b}\right)  ^{1/\left(  p-1\right)  }\delta_{\Omega
}\qquad\text{in }\overline{\Omega},
\]
and from the weak comparison principle (see Lemma \ref{w}) the same estimate
is also true for $\mathcal{S}\left(  \lambda^{\ast}h\right)  $. Furthermore,
by the homogeneity of the differential operator $\mathcal{L}$ we get that
\[
\mathcal{S}\left(  h\right)  \geq\left(  \frac{\min\left(  c_{a},c_{b}\right)
}{\lambda^{\ast}}\right)  ^{1/\left(  p-1\right)  }\delta_{\Omega}%
\qquad\text{in }\overline{\Omega}%
\]
which in turn yields (\ref{lap2}), and this ends the proof of the lemma.
$\blacksquare$

\begin{remark}
\label{hopf} Let us note that in particular (ii) establishes the strong
maximum principle and Hopf's Lemma for the operator $\mathcal{L}$, even if $h$
changes sign in $\Omega$. Moreover, it provides explicit lower and upper
bounds for $\mathcal{S}\left(  h\right)  ^{\prime}\left(  a\right)  $ and
$\mathcal{S}\left(  h\right)  ^{\prime}\left(  b\right)  $ respectively, in
terms of $\Omega$, $p$ and $h$. $\blacksquare$
\end{remark}

\qquad

Let $f:\Omega\times\left(  0,\infty\right)  \rightarrow\mathbb{R}$ be a
Carath\'{e}odory function (that is, $f\left(  \cdot,\xi\right)  $ is
measurable for all $\xi\in\left(  0,\infty\right)  $ and $f\left(
x,\cdot\right)  $ is continuous for $a.e.$ $x\in\Omega$). We consider next
singular problems of the form%
\begin{equation}
\left\{
\begin{array}
[c]{ll}%
\mathcal{L}u=f\left(  x,u\right)  & \text{in }\Omega\\
u>0 & \text{in }\Omega\\
u=0 & \text{on }\partial\Omega
\end{array}
\right.  \label{sing}%
\end{equation}
in a suitable sense. We say that $v\in W_{loc}^{1,p}\left(  \Omega\right)
\cap C\left(  \overline{\Omega}\right)  $ is a \textit{subsolution} (in the
sense of distributions) of (\ref{sing}) if $v>0$ in $\Omega$, $v=0$ on
$\partial\Omega$, and
\[
\int_{a}^{b}\left\vert v^{\prime}\right\vert ^{p-2}v^{\prime}\phi^{\prime}%
\leq\int_{a}^{b}f\left(  x,v\right)  \phi\qquad\text{for all }0\leq\phi\in
C_{c}^{\infty}\left(  \Omega\right)  .
\]
Analogously, $w\in W_{loc}^{1,p}\left(  \Omega\right)  \cap C\left(
\overline{\Omega}\right)  $ is a \textit{supersolution} of (\ref{sing}) if
$w>0$ in $\Omega$, $w=0$ on $\partial\Omega$, and
\[
\int_{a}^{b}\left\vert w^{\prime}\right\vert ^{p-2}w^{\prime}\phi^{\prime}%
\geq\int_{a}^{b}f\left(  x,w\right)  \phi\qquad\text{for all }0\leq\phi\in
C_{c}^{\infty}\left(  \Omega\right)  .
\]
For the sake of completeness we state the following existence theorem in the
presence of well-ordered sub and supersolutions (for the proof, see
\cite{loc}, Theorem 4.1).

\begin{theorem}
\label{subsup}Assume there exist $v,w\in C^{1}\left(  \Omega\right)  $ sub and
supersolutions respectively of (\ref{sing}), satisfying that $v\leq w$ in
$\Omega$. Suppose also that there exists $g\in L_{loc}^{p^{\prime}}\left(
\Omega\right)  $ such that $\left\vert f\left(  x,\xi\right)  \right\vert \leq
g\left(  x\right)  $ for $a.e.$ $x\in\Omega$ and all $\xi\in\left[  v\left(
x\right)  ,w\left(  x\right)  \right]  $. Then there exists $u\in C^{1}\left(
\Omega\right)  \cap C\left(  \overline{\Omega}\right)  $ solution (in the
sense of distributions) of (\ref{sing}) with $v\leq u\leq w$, that is,
\[
\int_{a}^{b}\left\vert u^{\prime}\right\vert ^{p-2}u^{\prime}\phi^{\prime
}=\int_{a}^{b}f\left(  x,u\right)  \phi\qquad\text{for all }\phi\in
C_{c}^{\infty}\left(  \Omega\right)  .
\]

\end{theorem}

\begin{remark}
\label{supersol}If $m\in L^{q}(\Omega)$ with $q>1$ and $m^{+}\not \equiv 0$,
one can quickly verify that (\ref{prob}) possesses arbitrarily big
supersolutions. Indeed, let $\psi:=\mathcal{S}\left(  m^{+}\right)  $ and let
us choose $\beta\in\left(  0,1\right)  $ and $\sigma>0$ satisfying%
\[
\beta:=\frac{p-1}{p-1+\gamma},\qquad\sigma\geq\frac{1}{\beta^{\beta}}.
\]
Notice that $\psi^{\beta}\in C^{1}\left(  \Omega\right)  \cap C\left(
\overline{\Omega}\right)  $, $\psi^{\beta}=0$ on $\partial\Omega$ and
$\psi^{\beta}>0$ in $\Omega$ by the strong maximum principle. Also, a simple
computation shows that
\begin{gather*}
\mathcal{L}\left(  \sigma\psi^{\beta}\right)  =-\left(  \sigma\beta\right)
^{p-1}\left(  \left\vert \psi^{\prime}\right\vert ^{p-2}\psi^{\prime}%
\psi^{\left(  \beta-1\right)  \left(  p-1\right)  }\right)  ^{\prime}=\\
\left(  \sigma\beta\right)  ^{p-1}\left(  m^{+}\left(  x\right)  \psi^{\left(
\beta-1\right)  \left(  p-1\right)  }-\left(  \beta-1\right)  \left(
p-1\right)  \left\vert \psi^{\prime}\right\vert ^{p}\psi^{\left(
\beta-1\right)  \left(  p-1\right)  -1}\right)  \geq\\
\left(  \sigma\beta\right)  ^{p-1}m^{+}\left(  x\right)  \psi^{\left(
\beta-1\right)  \left(  p-1\right)  }\geq m^{+}\left(  x\right)  \left(
\sigma\psi^{\beta}\right)  ^{-\gamma}\geq\\
m\left(  x\right)  \left(  \sigma\psi^{\beta}\right)  ^{-\gamma}\qquad\text{in
}\Omega^{\prime}%
\end{gather*}
for all $\Omega^{\prime}\Subset\Omega$, and hence $\sigma\psi^{\beta}$ is a
supersolution of (\ref{prob}). $\blacksquare$
\end{remark}

\section{Main results}

We denote%
\[
P^{\circ}:=\text{interior of the positive cone of }C_{0}^{1}\left(
\overline{\Omega}\right)  ,
\]
(that is, the functions $v\in C^{1}\left(  \overline{\Omega}\right)  $ with
$v\left(  a\right)  =v\left(  b\right)  =0$, $v>0$ in $\Omega$, $v^{\prime
}\left(  a\right)  >0$ and $v^{\prime}\left(  b\right)  <0$) and for any
$I=\left(  x_{0},x_{1}\right)  \subseteq\Omega$ we shall write
\[
x_{I}:=\frac{x_{0}+x_{1}}{2},\qquad c_{I}:=\max\left(  x_{I}-a,b-x_{I}\right)
.
\]

\begin{theorem}
\label{nuevo}Let $m\in L^{p^{\prime}}\left(  \Omega\right)  $ and $\gamma>0$.
\newline(i) Suppose
\begin{equation}
\mathcal{S}\left(  m\right)  \in P^{\circ}. \label{nece}%
\end{equation}
Then there exists $\gamma_{0}>0$ such that the problem (\ref{prob}) has a
solution $u\in P^{\circ}$ for all $\gamma\in\left(  0,\gamma_{0}\right]
$.\newline(ii) Suppose $m^{-}\delta_{\Omega}^{-\gamma}\in L^{q}\left(
\Omega\right)  $ with $q>1$. If for some $I=\left(  x_{0},x_{1}\right)
\subseteq\Omega$ it holds that%
\begin{gather}
\frac{\left(  \inf_{I}m^{+}\right)  ^{p-1+\gamma}}{\left(  \int_{a}^{b}%
m^{+}\right)  ^{\gamma}}\geq c_{\gamma,p,\Omega,I}\max\left(  \int_{a}^{x_{0}%
}m^{-}\delta_{\Omega}^{-\gamma},\int_{x_{1}}^{b}m^{-}\delta_{\Omega}^{-\gamma
}\right)  ^{p-1},\qquad\text{where}\label{cond}\\
c_{\gamma,p,\Omega,I}:=\left(  \frac{p-1}{\gamma}\right)  ^{\gamma}\left(
\frac{p-1+\gamma}{p-1}\right)  ^{p-1+\gamma}\left(  \frac{b-a}{2}\right)
^{\gamma\left(  p-1\right)  }\left(  c_{I}^{p-1}\lambda_{1}\left(  I\right)
\right)  ^{p-1+\gamma},\nonumber
\end{gather}
then the problem (\ref{prob}) has a solution $u\in C^{1}(\Omega)\cap
C(\overline{\Omega})$, and $u\in P^{\circ}$ whenever $m^{+}\delta_{\Omega
}^{-\gamma}\in L^{r}\left(  \Omega\right)  $ with $r>1$.
\end{theorem}

\textit{Proof}. Since Remark \ref{supersol} provides arbitrarily large
supersolutions of (\ref{prob}), it suffices to find a subsolution. Let us
start proving (i). We first observe that (\ref{prob}) admits a solution for
$m$ if and only if it has one for $\tau m$ for any constant $\tau>0$, and
therefore we shall also assume without loss of generality that $\mathcal{S}%
\left(  m^{+}\right)  \leq1$ in $\Omega$.

Due to (\ref{nece}), we can fix $\varepsilon>0$ such that $\mathcal{S}\left(
m\right)  \geq2\varepsilon\delta_{\Omega}$ in $\Omega$. We also pick
$\gamma_{0}>0$ such that for every $\gamma\in\left(  0,\gamma_{0}\right]  $ it
holds that $m^{-}\delta_{\Omega}^{-\gamma}\in L^{r}\left(  \Omega\right)  $
with $r>1$. Since $\mathcal{S}:L^{r}(\Omega)\rightarrow C^{1}(\overline
{\Omega})$ is a continuous operator for any $r>1$ (see Remark \ref{uno}),
making $\gamma_{0}$ smaller if necessary, we obtain that for all such $\gamma$
it holds that
\begin{equation}
\mathcal{S}\left(  m^{+}-m^{-}\left(  \varepsilon\delta_{\Omega}\right)
^{-\gamma}\right)  \geq\varepsilon\delta_{\Omega}\qquad\text{in }\Omega.
\label{aux}%
\end{equation}
Define now the set
\[
\mathcal{C}:=\left\{  v\in C\left(  \overline{\Omega}\right)  :\varepsilon
\delta_{\Omega}\leq v\leq\mathcal{S}\left(  m^{+}\right)  \text{ in }%
\Omega\right\}  ,
\]
and for $v\in\mathcal{C}$ let $u:=\mathcal{S}\left(  m^{+}-m^{-}v^{-\gamma
}\right)  :=\mathcal{T}\left(  v\right)  $. Utilizing (\ref{aux}) and the weak
comparison principle we see that%
\[
\mathcal{S}\left(  m^{+}\right)  \geq\mathcal{S}\left(  m^{+}-m^{-}v^{-\gamma
}\right)  =u\geq\mathcal{S}\left(  m^{+}-m^{-}\left(  \varepsilon
\delta_{\Omega}\right)  ^{-\gamma}\right)  \geq\varepsilon\delta_{\Omega
}\qquad\text{in }\Omega
\]
and hence $u\in\mathcal{C}$. Furthermore, one can verify that $v\rightarrow
m^{+}-m^{-}v^{-\gamma}$ is continuous from $\mathcal{C}$ into $L^{r}\left(
\Omega\right)  $ for some $r>1$, and thus employing the compactness of the
solution operator $\mathcal{S}$ (cf. Remark \ref{uno}) we deduce that
$\mathcal{T}:\mathcal{C}\rightarrow\mathcal{C}$ is continuous and compact. It
follows from Schauder's fixed point theorem that there exists some
$v\in\mathcal{C}$ solution of%
\begin{equation}
\left\{
\begin{array}
[c]{ll}%
\mathcal{L}v=m^{+}\left(  x\right)  -m^{-}\left(  x\right)  v^{-\gamma} &
\text{in }\Omega\\
v=0 & \text{on }\partial\Omega.
\end{array}
\right.  \label{auxx}%
\end{equation}
Moreover, $v\in C^{1}\left(  \overline{\Omega}\right)  $ and, since $v\leq1$
(due to $v\leq\mathcal{S}\left(  m^{+}\right)  \leq1$), it follows from
(\ref{auxx}) that $v$ is a subsolution of (\ref{prob}). Therefore, recalling
Remark \ref{supersol} and Theorem \ref{subsup} we obtain some $u\in
C^{1}\left(  \Omega\right)  \cap C\left(  \overline{\Omega}\right)  $ solution
of (\ref{prob}). Finally, decreasing $\gamma_{0}$ if necessary so that
$m^{+}\delta_{\Omega}^{-\gamma}\in L^{r}\left(  \Omega\right)  $ with $r>1$,
by standard regularity arguments we get that $u\in C^{1}\left(  \overline
{\Omega}\right)  $, and also $u\in P^{\circ}$ in view of the fact that $u\geq
c\delta_{\Omega}$ for some $c>0$. This concludes the proof of (i).

In order to prove (ii) we proceed similarly. We shall prove (ii) for $\tau m$,
where
\[
\tau:=\left(  \frac{2}{b-a}\right)  ^{p-1}\left(  \int_{a}^{b}m^{+}\right)
^{-1}.
\]
Since $\delta_{\Omega}\leq\left(  b-a\right)  /2$ in $\Omega$, employing
(\ref{lap1}) one can check that $\mathcal{S}\left(  \tau m^{+}\right)  \leq1$
in $\Omega$. We shall also assume that
\begin{equation}
\max\left(  \int_{a}^{x_{0}}m^{-}\delta_{\Omega}^{-\gamma},\int_{x_{1}}%
^{b}m^{-}\delta_{\Omega}^{-\gamma}\right)  =\int_{a}^{x_{0}}m^{-}%
\delta_{\Omega}^{-\gamma} \label{max}%
\end{equation}
because the other case is completely analogous. We define next%
\begin{gather*}
c_{1}:=\frac{\inf_{I}m^{+}}{\lambda_{1}\left(  I\right)  c_{I}^{p-1}},\qquad
c_{2}:=\int_{a}^{x_{0}}m^{-}\delta_{\Omega}^{-\gamma},\\
r:=\left(  \frac{\tau c_{2}\gamma}{p-1}\right)  ^{1/\left(  p-1+\gamma\right)
},\qquad\mathcal{C}:=\left\{  v\in C\left(  \overline{\Omega}\right)
:r\delta_{\Omega}\leq v\leq\mathcal{S}\left(  \tau m^{+}\right)  \text{ in
}\Omega\right\}  .
\end{gather*}
(Let us mention that if (\ref{max}) is not valid then we set $c_{2}%
:=\int_{x_{1}}^{b}m^{-}\delta_{\Omega}^{-\gamma}$.) One can readily verify
that (\ref{cond}) implies that%
\begin{equation}
c_{1}^{p-1+\gamma}\geq\left(  \frac{p-1}{\tau\gamma}\right)  ^{\gamma}\left(
\frac{p-1+\gamma}{p-1}\right)  ^{p-1+\gamma}c_{2}^{p-1}. \label{bueno}%
\end{equation}
Taking into account this fact and the definition of $c_{I}$ we now observe
that%
\begin{gather*}
\lambda_{1}(I)\left(  x_{I}-a\right)  ^{p-1}\int_{a}^{x_{0}}m^{-}\left(
r\delta_{\Omega}\right)  ^{-\gamma}\leq\lambda_{1}(I)c_{I}^{p-1}r^{-\gamma
}\int_{a}^{x_{0}}m^{-}\delta_{\Omega}^{-\gamma}=\\
\lambda_{1}(I)c_{I}^{p-1}\left(  \left(  \frac{p-1}{\tau\gamma}\right)
^{\gamma}c_{2}^{p-1}\right)  ^{1/\left(  p-1+\gamma\right)  }\leq\\
\lambda_{1}(I)c_{I}^{p-1}c_{1}\frac{p-1}{p-1+\gamma}<\inf_{I}m^{+}%
\end{gather*}
and thus we may apply Lemma \ref{plap} (ii) with $m^{+}-m^{-}\left(
r\delta_{\Omega}\right)  ^{-\gamma}$ in place of $h$ (and so also with
$\tau\left(  m^{+}-m^{-}\left(  r\delta_{\Omega}\right)  ^{-\gamma}\right)  $).

Given any $v\in\mathcal{C}$, we next define $u:=\mathcal{S}\left(  \tau\left(
m^{+}-m^{-}v^{-\gamma}\right)  \right)  $. Recalling the above paragraph, from
Lemma \ref{plap} (ii) and again making use of (\ref{max}) and (\ref{bueno}),
after some computations we deduce that%
\begin{gather*}
\mathcal{S}\left(  \tau m^{+}\right)  \geq u\geq\mathcal{S}\left(  \tau\left(
m^{+}-m^{-}\left(  r\delta_{\Omega}\right)  ^{-\gamma}\right)  \right)  \geq\\
\left(  \tau\left(  c_{1}-c_{2}r^{-\gamma}\right)  \right)  ^{1/\left(
p-1\right)  }\delta_{\Omega}\geq r\delta_{\Omega}\text{\qquad in }\Omega
\end{gather*}
and therefore $v\in\mathcal{C}$. Now the proof of (ii) can be finished as in
(i), and this concludes the proof of the theorem. $\blacksquare$

\begin{remark}
\label{dos} (i) Let us notice that by Lemma \ref{plap} (ii), (\ref{nece}) is
true if for instance
\[
\inf_{I}m^{+}>\lambda_{1}\left(  I\right)  \max\left(  \left(  x_{I}-a\right)
^{p-1}\int_{a}^{x_{0}}m^{-},\left(  b-x_{I}\right)  ^{p-1}\int_{x_{1}}%
^{b}m^{-}\right)
\]
for some $I=\left(  x_{0},x_{1}\right)  \subset\Omega$.\newline(ii) We also
remark that several distinct conditions guarantee that $m^{-}\delta_{\Omega
}^{-\gamma}\in L^{q}\left(  \Omega\right)  $ for some $q>1$. Indeed, for
example, this occurs for all $\gamma\in\left(  0,1/p\right)  $, or more
generally if $m^{-}\in L^{q}\left(  \Omega\right)  $ with $q\geq p^{\prime}$
and $\gamma\in\left(  0,\left(  q-1\right)  /q\right)  $. Also, the same is
valid for every $\gamma>0$ when $m\geq0$ in the set $\left\{  x\in
\Omega:\delta_{\Omega}\left(  x\right)  <\varepsilon\right\}  $ for some
$\varepsilon>0$. Of course analogous statements hold for $m^{+}\delta_{\Omega
}^{-\gamma}$. $\blacksquare$
\end{remark}

\begin{theorem}
\label{necesaria}Suppose (\ref{prob}) has a solution $u\in C^{1}\left(
\overline{\Omega}\right)  $ such that $\varphi_{p}\left(  u^{\prime}\right)  $
is absolutely continuous. Then
\begin{align}
\mathcal{S}\left(  m\right)   &  >0\qquad\text{in }\Omega\qquad\text{and}%
\label{ne}\\
\int_{a}^{b}m  &  >0. \label{refe}%
\end{align}

\end{theorem}

\textit{Proof}. Let $u>0$ be a solution of (\ref{prob}) and fix%
\[
\beta:=\frac{p-1+\gamma}{p-1}.
\]
Let $0\leq\phi\in C_{c}^{\infty}\left(  \Omega\right)  $, and let
$\Omega^{\prime}$ be an open set such that \textit{supp} $\phi\subset
\Omega^{\prime}\Subset\Omega$. We have that
\begin{gather*}
\mathcal{L}\left(  u^{\beta}\right)  =-\beta^{p-1}\left(  \left\vert
u^{\prime}\right\vert ^{p-2}u^{\prime}u^{\left(  \beta-1\right)  \left(
p-1\right)  }\right)  ^{\prime}=\\
\beta^{p-1}\left(  m\left(  x\right)  u^{-\gamma}u^{\left(  \beta-1\right)
\left(  p-1\right)  }-\left\vert u^{\prime}\right\vert ^{p}\left(
\beta-1\right)  \left(  p-1\right)  u^{\left(  \beta-1\right)  \left(
p-1\right)  -1}\right)  \leq\\
\beta^{p-1}m\left(  x\right)  u^{-\gamma}u^{\left(  \beta-1\right)  \left(
p-1\right)  }=\beta^{p-1}m\left(  x\right)  \qquad\text{in }\Omega^{\prime}%
\end{gather*}
and hence, multiplying the above inequality by $\phi$, integrating over
$\Omega^{\prime}$ and using the integration by parts formula we see that%
\[
\int_{a}^{b}\left\vert \left(  u^{\beta}\right)  ^{\prime}\right\vert
^{p-2}\left(  u^{\beta}\right)  ^{\prime}\phi^{\prime}\leq\beta^{p-1}\int
_{a}^{b}m\left(  x\right)  \phi.
\]

On the other hand, let $0\leq v\in W_{0}^{1,p}\left(  \Omega\right)  $. It is
easy to check that there exists $\left\{  \phi_{j}\right\}  _{j\in\mathbb{N}%
}\subset C_{c}^{\infty}\left(  \Omega\right)  $ with $\phi_{j}\geq0$ in
$\Omega$ and such that $\phi_{j}\rightarrow v$ in $W^{1,p}\left(
\Omega\right)  $ (see e.g. \cite{chipot}, p. 50). Utilizing the last
inequality with $\phi_{j}$ in place of $\phi$ and passing to the limit we get
that $\mathcal{L}\left(  u^{\beta}\right)  \leq\beta^{p-1}m\left(  x\right)  $
in weak sense in $\Omega$ and so from the weak comparison principle we deduce
that $0<u^{\beta}\leq\beta\mathcal{S}\left(  m\right)  $ in $\Omega$ and this
ends the proof of (\ref{ne}).

Finally, we observe that multiplying (\ref{prob}) by $u^{\gamma}$ and
integrating by parts on $\left(  a+\varepsilon,b-\varepsilon\right)  $ with
$\varepsilon>0$ small, we get that
\[
\left(  \varphi_{p}\left(  u^{\prime}\right)  u^{\gamma}\right)  \left(
a+\varepsilon\right)  -\left(  \varphi_{p}\left(  u^{\prime}\right)
u^{\gamma}\right)  \left(  b-\varepsilon\right)  +\gamma\int_{a+\varepsilon
}^{b-\varepsilon}\left\vert u^{\prime}\right\vert ^{p}u^{\gamma-1}\leq
\int_{a+\varepsilon}^{b-\varepsilon}m
\]
and letting $\varepsilon\rightarrow0$ it is easy to deduce (\ref{refe}).
$\blacksquare$

\begin{remark}
\label{rito} Let us note that the conditions (\ref{ne}) and (\ref{refe}) are
not comparable. Indeed, suppose first $p=2$, and let $\Omega:=\left(
0,3\pi\right)  $ and $m\left(  x\right)  :=\sin x$. Then $m=\mathcal{S}\left(
m\right)  $ and $\int_{0}^{3\pi}m>0$, but $\mathcal{S}\left(  m\right)  <0$ in
$\left(  \pi,2\pi\right)  $.

On the other side, integrating (\ref{g}) (with $m$ in place of $h$) we get
that
\begin{equation}
\varphi_{p}\left(  \mathcal{S}\left(  m\right)  ^{\prime}\left(  a\right)
\right)  -\varphi_{p}\left(  \mathcal{S}\left(  m\right)  ^{\prime}\left(
b\right)  \right)  =\int_{a}^{b}m. \label{fercho}%
\end{equation}
It follows that we may have $\mathcal{S}\left(  m\right)  >0$ in $\Omega$ but
$\int_{a}^{b}m=0$. (Take for instance again $p=2$, $\Omega:=\left(
0,\pi\right)  $, $m\left(  x\right)  :=2\left(  \sin^{2}x-\cos^{2}x\right)  $
and $\mathcal{S}\left(  m\right)  \left(  x\right)  =\sin^{2}x.$)

What it is indeed true from (\ref{fercho}) is that $\mathcal{S}\left(
m\right)  >0$ in $\Omega$ implies $\int_{a}^{b}m\geq0$. Moreover, from Theorem
\ref{necesaria} and (\ref{fercho}) we have that if (\ref{prob}) admits a
solution, then either $\mathcal{S}\left(  m\right)  ^{\prime}\left(  a\right)
\not =0$ or $\mathcal{S}\left(  m\right)  ^{\prime}\left(  b\right)  \not =0$.
It is an interesting open question to see if it is necessary that both
derivatives are nonzero. $\blacksquare$
\end{remark}

We conclude the paper showing an existence theorem for singular problems of
the form
\begin{equation}
\left\{
\begin{array}
[c]{ll}%
\mathcal{L}u=m\left(  x\right)  f\left(  u\right)  & \text{in }\Omega\\
u>0 & \text{in }\Omega\\
u=0 & \text{on }\partial\Omega,
\end{array}
\right.  \label{f}%
\end{equation}
for certain continuous functions $f:\left(  0,\infty\right)  \rightarrow
\left(  0,\infty\right)  $. Let us observe that we make no monotonicity nor
convexity assumptions on $f$.

We state the hypothesis

(H) There exist $c_{f},C_{f}>0$ and $\gamma>0$ such that%
\[
c_{f}\xi^{-\gamma}\leq f\left(  \xi\right)  \leq C_{f}\xi^{-\gamma}\text{ for
all }\xi>0.
\]

\begin{corollary}
\label{fff}Let $m\in L^{p^{\prime}}\left(  \Omega\right)  $, let $f$ satisfy
(H) and suppose (\ref{prob}) has a solution with $c_{f}m^{+}-C_{f}m^{-}$ in
place of $m$. Then there exists a solution of (\ref{f}).
\end{corollary}

\textit{Proof}. Let $u\ $be a solution of (\ref{prob}) with $c_{f}m^{+}%
-C_{f}m^{-}$ in place of $m$. Employing (H) we find that
\[
\mathcal{L}u=\left(  c_{f}m^{+}\left(  x\right)  -C_{f}m^{-}\left(  x\right)
\right)  u^{-\gamma}\leq m\left(  x\right)  f\left(  u\right)  \qquad\text{in
}\Omega\text{.}%
\]
On the other hand, let $\psi:=\mathcal{S}\left(  m^{+}\right)  >0$ and fix
$\beta\in\left(  0,1\right)  $ and $\sigma>0$ satisfying%
\[
\beta:=\frac{p-1}{p-1+\gamma},\qquad\sigma\geq\frac{C_{f}^{1/\left(
p-1+\gamma\right)  }}{\beta^{\beta}}.
\]
Enlarging $\sigma$ if necessary, recalling that $\beta<1$ and that by Lemma
\ref{plap} $\mathcal{S}\left(  m^{+}\right)  ^{\prime}\left(  a\right)
>0>\mathcal{S}\left(  m^{+}\right)  ^{\prime}\left(  b\right)  $, we may
assume that $\sigma\psi^{\beta}\geq u$ in $\Omega$. Now, arguing as in Remark
\ref{supersol} and taking into account (H) we obtain that%
\begin{gather*}
\mathcal{L}\left(  \sigma\psi^{\beta}\right)  \geq\left(  \sigma\beta\right)
^{p-1}m^{+}\left(  x\right)  \psi^{\left(  \beta-1\right)  \left(  p-1\right)
}\geq\\
C_{f}m^{+}\left(  x\right)  \left(  \sigma\psi^{\beta}\right)  ^{-\gamma}\geq
m^{+}\left(  x\right)  f\left(  \sigma\psi^{\beta}\right)  \geq m\left(
x\right)  f\left(  \sigma\psi^{\beta}\right)  \qquad\text{in }\Omega^{\prime}%
\end{gather*}
for every $\Omega^{\prime}\Subset\Omega$, and the corollary follows.
$\blacksquare$

\end{document}